\documentclass{article}\begin{document}\centerline{\bf Meromorphic Reduction in Integration}\vskip 1in

\centerline{ Lawrence Glasser}
\centerline{ Professor Emeritus}

\centerline{Department of Physics, Clarkson University}

\centerline{Potsdam, NY 13699-5820}  
\centerline{e-mail: laryg@tds.net}
\vskip 1.5in   
\centerline{\bf ABSTRACT}
\begin{quote}

It is argued that for certain meromorphic functions $u:\cal{R}\rightarrow\cal{R}$ and analytic function $ A_1$ and  for any integrable function $F$, as long as it converges as a Cauchy Principal Value,,

$$\int_{-\infty}^{\infty}A_1(x)F[u(x)] dx=\int_{-\infty}^{\infty} A_2(x)F(x) dx,$$
where $A_2$ is also analytic.

\end{quote}\vskip ,4in
{\bf Keywords: Meromorphic function, analytic function, definite integral.\vskip .5in
\noindent
{\bf MSC:  26A06,  26A42}

\newpage

\centerline{{\bf 1. Introduction}}
\vskip .1in
About 40 years ago the author published a generalization of the famous Cauchy-Schl\"omilch identity [1],[2], that for any integrable real valued function $F$,
$$\int_{-\infty}^{\infty}F\left(x-\frac{1}{x}\right)dx=\int_{-\infty}^{\infty}F(x)dx\eqno(1)$$
subject only to the minor restriction, that if there are simple poles in the integrand, the Cauchy principal value prescription is appropriate. It will be shown in the next section that, subject only to integrability, for any real sequences    $\{a_k\},\;\{b_k\}$, $N\in {\cal Z}$ and polynomial $p(x)$ 
$$\int_{-\infty}^{\infty} p(x)F\left(x-\sum_{k=1}^N\frac{|a_k|}{x-b_k}\right)dx=\int_{-\infty}^{\infty}q(x)F(x)dx,\eqno(2)$$
where $q$ is a polynomial of the same degree as $p$. Now, as long as the $p$ and $q$   polynomials converge and the integrands remain integrable as the degree of $p$ becomes infinite, they are, {\it de facto} analytic functions. Assertion (2) need only be proven for $p(x)=x^n, n\in \cal{Z}^+$ as will be carried out explicitly in Section 2 for $n=1,2,3,4$; it should become clear that there is no barrier to larger values. Several examples will be worked out in Section 3.\vskip .3in

\centerline{\bf 2. Calculation}\vskip .2in

Let $\{a_k\}$ be any sequence of positive real numbers, $\{b_k\}$ be any sequence of reals (having the same lengths $n$) and let

$$u=x-\sum_{k=1}^n\frac{a_ k}{x-b_k}\eqno(3)$$ 
A plot of u vs x consists of $n+1$  continuous monotonic segments  separated by the  poles $b_1,\cdots, b_n$. Each segment has a unique differentiable  inverse $x_k(u)$.
 Hence, for any integrable function $F$
            $$I_m=\int_{-\infty}^{\infty}x^m F(u)dx=\int_{-\infty}^{b_1}x^mF(u)dx+\int_{b_1}^{b_2}x^mF(u)dx+\cdots +\int_{b_n}^{\infty}x^mF(u)dx.\eqno(4)$$
      Now, by changing the integration variable to u in each section, this becomes
      $$I_m=\int_{-\infty}^{\infty}\sum_{k=1}^n x_k^m x_k'(u)F(u)du\eqno(5)$$
      where the derivative is with respect to $u$.
      
      Next, from (3) the equation for the functions $x_j(u)$
     $$ (u-x)\Pi_{k=1}^n(x-b_k)+\sum_{k=1}^n a_k{\bar \Pi}_{l=1}^n(x-b_l)=x^{n+1}+\sum_{k=0}^n(-1)^k\sigma_{k+1}x^{n-k}=0\eqno(6)$$
     where the barred product means $l\ne k$ and the $\sigma_k$ denote the elementary symmetric polynomials in the $x_k$.  We  require the power sums
     $$\tau_j=\sum_{k=1}^nx_k^j, \quad j=1,2,3,\cdots\eqno(7)$$
     which are related to the $\sigma_j$ by Newton's identities[3]. ( Actually, Albert Girard,  1595-1632, who also  introduced  the notation {\it sin, cos and tan}, found these relations many years before Newton).
      
      $$\tau_1=\sigma_1,\quad \tau_2=\sigma_1^2-2\sigma_2,\quad \tau_3=\sigma_1^3-3\sigma_1\sigma_2+3\sigma_3$$
      $$\tau_4=\sigma_1^4-4\sigma_1^2\sigma_2+4\sigma_1\sigma_3-2\sigma_2^2 +4\sigma_4\eqno(8)$$
      $$\tau_5=\sigma_1^5-3\sigma_1^3\sigma_2+5\sigma_1^2\sigma_3+\sigma_1\sigma_4-\sigma_1\sigma_2^2+\sigma_2\sigma_3+\sigma_1\sigma_4+5\sigma_5$$
      For  equation (6), one has by inspection
      $$\sigma_0=1,   \quad \sigma_1=\beta_1+u, \quad \sigma_2=\alpha_1-\beta_2-u\beta_1\eqno(9)$$
      (working out  the subsequent  relations by hand becomes increasingly tedious). where, $\alpha_j$ and $\beta_j$ are the symmetric polynomials for the  sequences $a_k$ and $b_k$, resp.
      
      Returning to (5), we now have
      $$  \int_{-\infty}^{\infty}x^m F(u)dx=  \int_{-\infty}^{\infty }       q_m(u)F(u)du,   \quad  q_m(u)=\frac{1}{m+1}\frac{d}{du}\tau_{m+1} \eqno(10)$$
      and $q_m(u)$, which   is of  degree $m$ in $u$, is explicitly expressible in terms of the $\sigma'$s which are polynomials in $u$.
        By making use of (8) and (9) one finds, for $m=0,1,2,3$
        
      $$  \int_{-\infty}^{\infty}F(u)dx=\int_{-\infty}^{\infty}F(x)dx\eqno(11)$$
      $$ \int_{-\infty}^{\infty}xF(u)dx=\int_{-\infty}^{\infty}xF(x)dx\eqno(12)$$
      $$ \int_{-\infty}^{\infty}x^2F(u)dx=\int_{-\infty}^{\infty} (x^2+\alpha_1)F(x)dx\eqno(13)$$
      $$ \int_{-\infty}^{\infty}x^3F(u)dx=\int_{-\infty}^{\infty} (x^3+2\alpha_1x+\sum_{k=1}^na_kb_k)F(x)dx\eqno(14)$$
      $$\cdots$$
      \vskip .3in
      
      \centerline{\bf 3. Examples}\vskip .2in
      
      As a first example, consider the untabulated integral
      
      $$I_1=\int_0^{\infty}x^{-5/2}\exp\left[-\frac{1}{x}\left(\frac{1-ax}{1-bx}\right)^2\right]dx,   \quad 0<b<a.\eqno(15)$$
      First, to eliminate the fractional power, replace $x$ by $v^2$. Then
      $$I_1=2\int_0^{\infty} v^{-4}\exp\left[-\left(\frac{1-av^2}{v(1-bv^2)}\right)^2\right]dv.\eqno(16)$$
      Next, to reduce the exponent in the integrand, let $v\rightarrow 1/t$
      
      $$I_1=2\int_0^{\infty} t^2\exp\left[-\left(t\frac{t^2-a}{t^2-b}\right)^2\right]dt  $$
   $$  = \int_{-\infty}^{\infty} t^2 \exp\left[-\left(t-\frac{(a-b)/2}{t-\sqrt{b}}-\frac{(a-b)/2}{t+\sqrt{b}}  \right)^2\right]dt\eqno(17).$$
   Finally, by applying (13) we get the known integral
   $$I_1=2\int_0^{\infty}(x^2+a-b)\exp(-x^2)dx=\sqrt{\pi}(a-b+1/2).\eqno(18)$$
   
   As a second example, consider
   $$I_2= \int_0^{\infty} x\;  {\rm csch}\left(x-\frac{\pi}{\tan(\pi x)}\right)dx\eqno(19)$$
   Since
   $$\frac{\pi }{\tan(\pi x)}=\frac{1}{x}+\sum_{k=1}^{\infty}\left(\frac{1}{x-k}+\frac{1}{x+k}\right),\eqno(20)$$
   by (12) and symmetry we have
  $$I_2= \int_0^{\infty} x\; {\rm csch}(x)dx=\frac{\pi^2}{4}\eqno(20)$$
  However, in this case, since $\alpha_1$ and $\beta_1$ diverge, (19) cannot be extended to higher powers of $x$.
  
  A final example is (where $a^2=2\pi+9$)
  $$I_3=\int_0^{\infty}( x^4+4x^2+1)\exp\left[-2x^2\left(\frac{x^2-a^2}{x^2-9}\right)^2\right].  \eqno(21)$$
  First note that  the exponential argument is $2u^2$ where  $u-x+\pi[(x-3)^{-1}+(x+3)^{-1}=0]$.  For this equation 
  $$\tau_2=u^2+4\pi+18,\quad \tau_3=u^3+6\pi u \; \mbox{ and  }\;  \tau_5=u[u^4+10\pi(u^2+9+20\pi^2].\eqno(22)$$
  $$q_4(u)=u^4+6\pi u^2+4\pi^2+18\pi),\quad q_2(u)=u^2+2\pi \mbox{. and } q_0(u)=1\eqno(23)$$.
   Hence
  $$I_3=\int_0^{\infty}[q_4(u)+4q_2(u)+1)e^{-2u^2} du=\sqrt{\frac{\pi}{8}}\left(\frac{35}{16}+4\pi^2+\frac{11}{4}\right).\eqno(24)$$
  \vskip .3in
  
  \centerline{\bf Conclusion}\vskip .1in
  
  A search of the most comprehensive contemporary compilations of integral identities identified only five entries related to the present work. These are all listed together in section 2.1.2, volume 1 of the five volume set of integral tables [4] as entries 49-53. Entry 49 is incorrect  and they all are for the Cauchy- Schl\"omilch case $u=ax-b/x$.  As pointed out in [1] meromorphic reduction is also applicable for 
  $$u_1=\sum_{k=1}^n\frac{a_k}{x-b_k} \eqno(25)$$
  It  is important to note that in 1977 G.  Letac[5] pointed out that  meromorphic reduction can be related to the {\it  Cauchy Law}  in statistics   (but, only for $p=1$) if $u$ is measure preserving. Since the iteration of a measure preserving transformation is measure preserving, we find such bizarre  sequences as
  $$\sqrt{\pi}        =\int_{-\infty}^{\infty}\exp{\left[-\left( \frac{x^2-1}{x}  \right)^2\right]}dx   =\int_{-\infty}^{\infty}
  \exp{\left[-\left(\frac{1-2x}{x(x^2-1)}\right)^2\right]} dx =\cdots     \eqno(26) $$
  It is an open question whether the same is true for meromorphic reduction. Another open question is whether these identities possess higher dimensional analogues.

      \centerline{\bf References}\vskip .1in
      \noindent 
      [1] M.L.  Glasser, {\it A Remarkable Property of Definite Integrals}, Math. Comp.{\bf 40}, 561-563 (1983)
      
      \noindent
      [2] T. Andeberhan, M.L. Glasser, M.C. Jones, V. Moll, R.Posey  and  D. Varela,{\it  The Cauchy-Schl\"omilch Transformation. },Int. Trans. and Special Functs.
       {\bf 30}  940-961 (2019).
       
       \noindent
       [3] B.L. Van der Waerden,{\it Modern Algebra}[ Frederick Ungar Publishing Co, New York,(1931) ], Vol 1, Ch.4

\noindent
[4] A.P. Prudnikov, Yu.A. Brychkov and O.I. Marichev,{\it Integrals and Series},[Gordon and Breach Science Publishers, NEW YORK (1986)]

\noindent
[5]   G. Letac,{\it  Which Functions Preserve Cauchy's Law?}, Proc. Amer. Math. Soc.  {\bf 67}, 277-286 (1977).

\end{document}